\pdfoutput=1 
\documentclass[12pt]{article}
\usepackage{amsmath,amsfonts,amssymb,url,color,epsfig,graphics}
\usepackage[latin1]{inputenc}

\newcommand{\R}{{\mathbb{R}}}

\newcommand{\Z}{{\mathbb{Z}}}

\newcommand{\D}{{\mathbb{D}}}
\newcommand{\T}{{\mathbb{T}}}
\newcommand{\F}{{\mathbb{F}}}

\def\ha{\frac{1}{2}}

\def\pa{\partial}
\def\ra{\rightarrow}
\def\ga{\alpha}

\def\ge{\varepsilon}

\def\gg{\gamma}

\def\gl{\lambda}
\def\go{\omega}

\def\gs{\sigma}

\def\nghbd{neighbourhood~}
\def\san{San V{\~u} Ng{\d o}c}
\newtheorem{defi}{Definition}[section]

\newtheorem{lemm}{Lemma}[section]
\newtheorem{prop}{Proposition}[section]

\newtheorem{coro}{Corollary}[section]
\newtheorem{theo}{Theorem}[section]

\begin{document}

\title{Large time asymptotics of    the  wave fronts length II:\\
surfaces with  integrable Hamiltonians}

\author{Yves  Colin de Verdi\`ere\footnote{Universit\'e Grenoble-Alpes,
Institut Fourier,
 Unit{\'e} mixte
 de recherche CNRS-UGA 5582,
 BP 74, 38402-Saint Martin d'H\`eres Cedex (France);
{\color{blue} {\tt yves.colin-de-verdiere@univ-grenoble-alpes.fr}}}}

\maketitle
In the paper \cite{Vi-20}, the  author proves that the length $|S_t|$ of the wave front $S_t$ at time $t$ of a wave propagating in an
Euclidean  disk $\D $ of radius $1$,
starting from a source   $A$, admits a linear asymptotics as $t\ra + \infty $:
$|S_t|\sim (2 \arcsin a) t $   with $a=d(0,A)$. 
In the paper \cite{Co-Vi-20}, we gave a more direct proof and some improvements of that
result.

Here, we will  explain that this result is quite general for  surfaces 
with an {\it integrable Hamiltonian}. 
We  discuss only   the 2D case for simplicity. 
The main idea  is to use {\it action-angle coordinates} (section \ref{sec:morse-bott})
in order to  get a  nice integral expression for  $|S_t|$ (section \ref{sec:A}).
Integrable systems have in  general singularities, therefore we need to make some genericity assumptions
(section \ref{sec:morse-bott})
and to study  what happens to the action-angle coordinates (section \ref{sec:singu})
near these generic  singularities. We need then to  evaluate  some
oscillatory integrals (section \ref{sec:proof}) using an ergodic lemma (Appendix \ref{app:ergo}).

For the geodesic flow on closed 
manifolds of negative curvature, Margulis  \cite{Ma-69} proved that the asymptotics of the length is exponential.
The generic behaviour is not known. Here we study the integrable case which
is highly non generic. 

Before starting, let us give a rough  version of the main theorem \ref{theo:main2}:

{\bf Let $(X,g)$ be a 2D-Riemannian manifold. Let $H:T^\star X \ra \R $ be an 
  integrable Hamiltonian near a given  energy $E$. Assume that the  energy shell $\Sigma:=H^{-1}(E) $ is compact
and 
that $dH $ does not vanish on $\Sigma $.
We assume also that $H$ satisfies some ``generic properties''.    The wave front starting from a point $A$ at energy $E$
is the projection onto $X$ of $\phi_t (\Sigma^A )$ 
where $\phi_t$ is the Hamiltonian flow of $H$ 
and $\Sigma^A:=\{(A,\xi)|H(A,\xi)=E \}$ is   assumed to be smooth.  If the point $A$  in $X$ is ``generic'' 
(see section \ref{sec:A}),  the g-length of the wave front starting from $A$, at energy $E$, admits a linear asymptotics
$|S_t| \sim \lambda (A) t $ as $t \ra +\infty $, where $\lambda (A)$  expresses in terms of the action-angle coordinates.}


\section{Wave fronts}
Let us consider a smooth 2D  Riemannian   manifold $(X,g)$  without boundary and fix a real number $E$. Let
$H:T^\star X  \ra \R $ be a smooth  Hamiltonian. Assume that $H^{-1}([E-a, E+a ])$ is compact for
some positive $a$ and that $dH $ does not vanish on 
 $\Sigma :=H^{-1}(E)$. Let us fix some point $A \in X$ and put $\Sigma^A:=\{ (A,\xi)\in \Sigma \}$. We denote by $\go $
the generic point of $\Sigma^A$.  Assume that 
$d(H_{|T^\star_A X}) $ does not vanish on $\Sigma^A $. This implies that $\Sigma^A$ is a 1D-compact submanifold of $\Sigma $. 
 We denote by $\pi _X$ the canonical projection
of $\Sigma   $ onto $X$ and by $\phi_t:\Sigma  \ra \Sigma ,~t\in \R$, the   flow of $\vec{H}$, the Hamiltonian vector field
derived from $H$.   
 For any positive $t$, we define the {\it wave front} $S_t$ at time  $t$  as the set of points of $X$ 
of the form $\pi_X (\phi_t (\Sigma ^A ) )$. The wave front $S_t$ has a smooth  parametrization by $\Sigma ^A$.
This allows to define its length $|S_t|$ using the Riemannian metric $g$, assumed to be continuous and possibly degenerate:
\[ |S_t|=\int_{\Sigma^A}\gg^{\ha}\left(\phi_t(\go); \frac{d}{d\go}\phi_t(\go)\right)|d\go| \]
where $\gg=\pi_X^\star (g)$. 
 Note that $S_t$  admits in general some singular points as a subset of $X$,
namely cusps and transversal self-intersections. 
We are interested 
in the asymptotic behaviour of $|S_t|$ as $t\ra +\infty $. 

{\it Examples:}
\begin{enumerate}
\item
 {\it Geodesic flows:} $H:=\ha g^\star $  is the Hamiltonian of the geodesic flow of a  closed Riemannian manifold
$(X,g)$. 
Let us fix $E=2$. Then $\Sigma $ is the unit cotangent bundle and, on $\Sigma$,  $\phi_t$ is the geodesic flow with speed $1$.  
In this case, $S_t$ is the image by the exponential map at the point $A$ of the circle of radius $t$ in the tangent plane
$T_AX$. 
\item
 {\it Schr\"odinger Hamiltonians:} Let $(X,g)$ be a  Riemannian manifold without boundary, 
$V:X \ra \R $ a smooth function and $E$ a real number.
 We take $H:=\ha g^\star + V $. 
Our  assumptions are satisfied  if
 $V^{-1}(]-\infty ,E+a])$ is compact for some $a>0$, 
  $dV $ does not vanish on $V^{-1}(E)$ and 
$V(A)<  E$. \end{enumerate}
\section{Integrable Hamiltonian  flows}\label{sec:morse-bott}
For this section, one can look at the chapter 4 of \cite{Vu-06} and the section 1 of  \cite{Co-Vu-03}.
We will assume  that the Hamiltonian $H$   is  integrable near the energy $E$. ``Integrability''   means
 that there exists a positive number $a$ and 
a smooth map $M=(I,J):H^{-1}(]E-a, E+a[)\ra \R^2 $, called the {\it moment map}, so that
\begin{itemize}\item The Poisson bracket $\{ I,J\} $ vanishes identically.
\item The critical points of $M$ are of measure $0$, i.e. the differentials  $dI $ and $dJ$ are almost everywhere independent.
\item There exists a smooth function $\Phi: M(H^{-1}(]E-a, E+a[))\ra \R $ so that
$H=\Phi(I,J)$. 
\end{itemize}
Note that $dI$ and $dJ$ cannot vanish at the same point of $\Sigma $ because $dH$ does not vanish there. 

The main examples with the geodesic flows are the surfaces of revolution, the tri-axial ellipso\"ids (\cite{Ja-39})
and the Liouville metrics on 2D tori
(Liouville metrics are of the form $ds^2=(f(u) +g(v))(du^2+ dv^2)$, see \cite{B-S-K-97} Chap. 7).
Usually, integrable systems have singularities. We will make the following ``generic'' assumption which is already used 
in \cite{Co-Vu-03}: we assume that the moment map  $M$ satisfies the 

 {\bf (A1) Morse-Bott condition: at any point of $\Sigma $ where $dI$ and $dJ$ are linearly dependent, i.e.
where $\gl dI+\mu dJ=0 $ for some  pair $(\gl,\mu)\ne (0,0)$, the function $\gl I +\mu J $, restricted to $\Sigma $, 
 admits a critical manifold of dimension $1$ with a transversally
non degenerate Hessian.}

 This implies that the singular set $Z_0\subset \Sigma $, i.e.
 the set of  critical points  of $M$  located in  $\Sigma $, is a finite
union of periodic orbits   of $\vec{H}$. These periodic orbits  are either hyperbolic or elliptic according to the signature of
the transversal Hessian.
We denote by $Z\subset \Sigma $ the part of the preimage by  $M$ of the critical values of $M$
which is the union  of $Z_0$ and all   the stable
and unstable manifolds of the hyperbolic  periodic orbits. The  open set $\Sigma \setminus Z_0$ admits  a smooth Lagrangian foliation given
by the level sets of $M$.

The open set $\Sigma \setminus Z $ is foliated by 2D-tori on which the Hamiltonian  flow of $H$ is quasi-periodic.
The set of these tori is a smooth 1D-manifold. We denote it  by ${\cal L}$ and  by $\gs $ the generic point  of ${\cal L}$.
 The manifold ${\cal L}$ is a 1D-torus in the case where there are no singularities, i.e. if $Z$ is empty,  and
a finite  union of real lines $D_j,~j=1,\cdots N$ if there are some singularities. 
 If $\gs \in D_j $ tends to one of the infinity of $D_j$, the corresponding torus $\T _\gs$
converges  to a  compact connected set  $\T_{j, \infty }$  of $\Sigma $ which is either  an elliptic periodic orbit of $H$ or
 the union of a finite set  of hyperbolic periodic orbits
 of $\vec{H}$ and  some  cylinders which are  connected components of their stable manifolds. In the last case,
$\T_{j, \infty }$ is homeomorphic to a 2D torus or to a Klein bottle.  

Let us denote by  $U_j $ the open connected component of
$\Sigma \setminus Z$ which is the union of the tori associated to the line $D_j$.
 The projection of $U_j$ onto $D_j$ is a  smooth fibration by
2D-tori which is trivial, because it is a fibration on the real line.
 There exist global coordinates $(\theta,\gs)\in \T^2 \times {\cal L}$ on $\Sigma \setminus Z$ so that the torus
$\T_\gs $ is mapped onto $\T^2 \times \{ \gs \}$ and the Hamiltonian flow is mapped on a 
vector field $V(\gs)=A(\gs)\pa _{\theta _1}+ B(\gs )\pa_{\theta _2} $ on $\T^2$
with some smooth functions $A$ and $B$. Note that $A$ and $B$ have no common zeroes because the Hamiltonian flow  does not vanish
on $\Sigma $.  

In what follows, we fix some component $U_j$. Let us describe the {\it action-angle coordinates} in some \nghbd  of  $U_j$
in $T^\star X$:
  there exists a    symplectic diffeomorphism $\chi_j $ of some  \nghbd  $V_j$ of $U_j$
   onto an open set $\T^2 \times \Omega_j$, with $\Omega _j \subset \R^2$,  contained in  
   $T^\star \T^2\setminus 0$ with canonical coordinates $(\theta,p)$,  so that 
$H \circ  \chi_j ^{-1}(\theta ,p) =K_j (p)$ with $K_j$ a smooth function  from $\Omega_j$ into $\R$. 
In these coordinates, the vector field $V_j$ is given
by $V_j=(\pa K_j /\pa p_1) \pa _ {\theta_1} +(\pa K_j /\pa p_2) \pa _ {\theta_2} $.
We note $\tilde{\nabla} K$ this non vanishing vector field.  The vector field  $\tilde{\nabla} K_j$ does not vanish and hence the curve 
 $C_j:=\{p \in \R^2|K_j(p)=1\}$ is a smooth submanifold of $\Omega _j $. 
The line $D_j$ identifies smoothly to the  curve $C_j$.
The manifold ${\cal L}$ can be identified to  the disjoint union of the curves $C_j$. 
The coordinates $p$ are called the actions: they are given by action integrals $p_j:=\int_{\gg _j}\ga $, where $d\ga $ is the symplectic 
form,
and the loops $\gg_j,~j=1,2$ form  a basis of $H_1(\T_\gs,\Z)$ varying continuously in $V_j$. Note that if $\ga'$ is another primitive
of the symplectic form, the difference 
$\ga -\ga'$ is closed, hence the action integrals differ by some constants. 
There are  many   choices for the coordinates $\theta $: if $\Lambda \subset U_j$  is a  Lagrangian manifold transversal
to the foliation by the tori, one can choose  $\theta $ vanishing  on $\Lambda $.

We will need one more ``generic'' assumption on the Hamiltonian flow:

 {\bf (A2) For any $j=1,\cdots, N$, 
 there exists, at any point $p$ of $C_j$, two integers $k\geq 1$ and $l\geq 1$, 
so that the derivatives of order $k$ and $l$ of the vector field ${\nabla }K_j $ along $C_j$ are linearly independent.}

Note that this condition is independent of the parametrization of $C_j$. 
For example, $k=1,l=2$ means that the curvature of the curve $\{ \nabla K _j(p)|p\in C_j \}$ does not vanish while
$k=2,l=3$ means a generic cusp for that curve. 

The assumption {\bf (A2)} implies that, for any $\nu \in \R^2 \setminus 0$,  the map from $C_j$ into $\R$
defined by $p \ra \langle \nu |\nabla K_j \rangle $ has only critical points of finite order.

\section{The behaviour of ${\nabla }K $ near $Z$}\label{sec:singu}
In this section, we forget about the index $j$: $K$ denotes the expression of $H$ in some of the action-angle coordinates.
 We are interested at the behaviour of $\nabla K$
near $ Z$. 
\subsection{The Elliptic case}

\begin{lemm} \label{lemm:gradK-ell}
Let $\gg $ be an   elliptic periodic orbit of $H$ (included in $Z_0$), then 
 $K$ is a smooth  function of  $p_1$ and $p_2$ up to $\gg$.
\end{lemm}
{\it Proof.--}
There exists a symplectic chart of a neighborhood of the elliptic periodic orbit of $H$ 
so that  $H=\Phi \left(\xi, y^2+\eta ^2 \right)$
in $(T^\star \T)_{x,\xi} \times (T^\star \R )_{y,\eta}$ (see \cite{Vu-00}).
The invariant tori are the level surfaces of the moment function $M(x,\xi;y,\eta)=(\xi,y^2 +\eta^2)$.  
Let us choose  $\gg_1= \{ s\ra (s,\xi;y,\eta )|s\in \R/\Z \}$ and
 $\gamma _2=\{ s \ra (x,\xi, \sqrt{y^2+\eta^2}\cos 2\pi s, \sqrt{y^2+\eta^2}\sin 2\pi s)|s\in\R/\Z \}$. If
$\ga =\xi dx +\eta dy$,  we get
the action integrals $p_1=\xi$ and $p_2= \pi (y^2+\eta^2) $. Hence $K(p_1,p_2) \equiv \Phi (p_1, p_2/\pi ) $.
\hfill $\square$
 
We have the following  
\begin{coro}\label{cor:int-nablaK-ell}
The manifold ${\cal L} $ admits an extension as a manifold with boundary at the elliptic periodic orbits of $H$ and 
the 1-form $d{\nabla }K  $ is smooth  and hence integrable on  ${\cal L}$ near  that boundary.
\end{coro}

\subsection{The Hyperbolic case}
In this section, we will use Section 1 of \cite{Co-Vu-03}.
\subsubsection{Functions of type {\bf (L)}}

Let us start with a
\begin{defi} A function $f:[0,c[ \ra \R$,  with $c>0$, is called of type {\bf (L)} if there exists two  smooth functions
$\phi, \psi :[0,c [ \ra \R $ so that 
\[\forall x\in [0,c[,~ f(x)=\phi(x) \log x + \psi(x) \]
with $\phi (0)=0,~\phi'(0)\ne 0 $.
\end{defi}
This definition is invariant by any smooth change of variable from $[0,c[$ into $[0,c'[$. Hence it extends to
$1D$-manifolds with boundaries. 
Such a function is  invertible in a small enough subinterval of $[0,c[$  and the inverse
$f^{-1}$ is $C^1$ up to the boundary $\psi(0)=\lim _{x\ra 0^+}f(x)$.  

Now let us describe the application that we have in mind.
\begin{lemm}\label{lemm:typel}
Let $\ga $ be a smooth 1-form so that $d\ga $ is a volume form in some  \nghbd $[-d,d]^2 $ of  the origin in the $(y,\eta)$ plane. 
Let $m_j(s)=(y_j(s), \eta_js)),~j=1,2$,  be two smooth curves with $m_1(0)=(e,0)$, 
$m_2(0) =(0,f) $ and $ e,f \in ]0,d[ $, $\eta'_1(0)>0$, $y_2'(0)>0 $ which are arcs transverse to each of the coordinates axes.
 Then consider
the integral $I(t)=\int _{\Gamma _t } \ga $ where $\Gamma _t $ is, for $t$ small enough, the part of the hyperbola $y\eta =t$ ($t>0$)
between the curves $m_1 $ and $m_2$ oriented in any of the two possible  directions. 
Then $I(t)$ is a function of type {\bf (L)}.\end{lemm}
\begin{figure}[!ht]
\centering
\includegraphics[height=7cm]{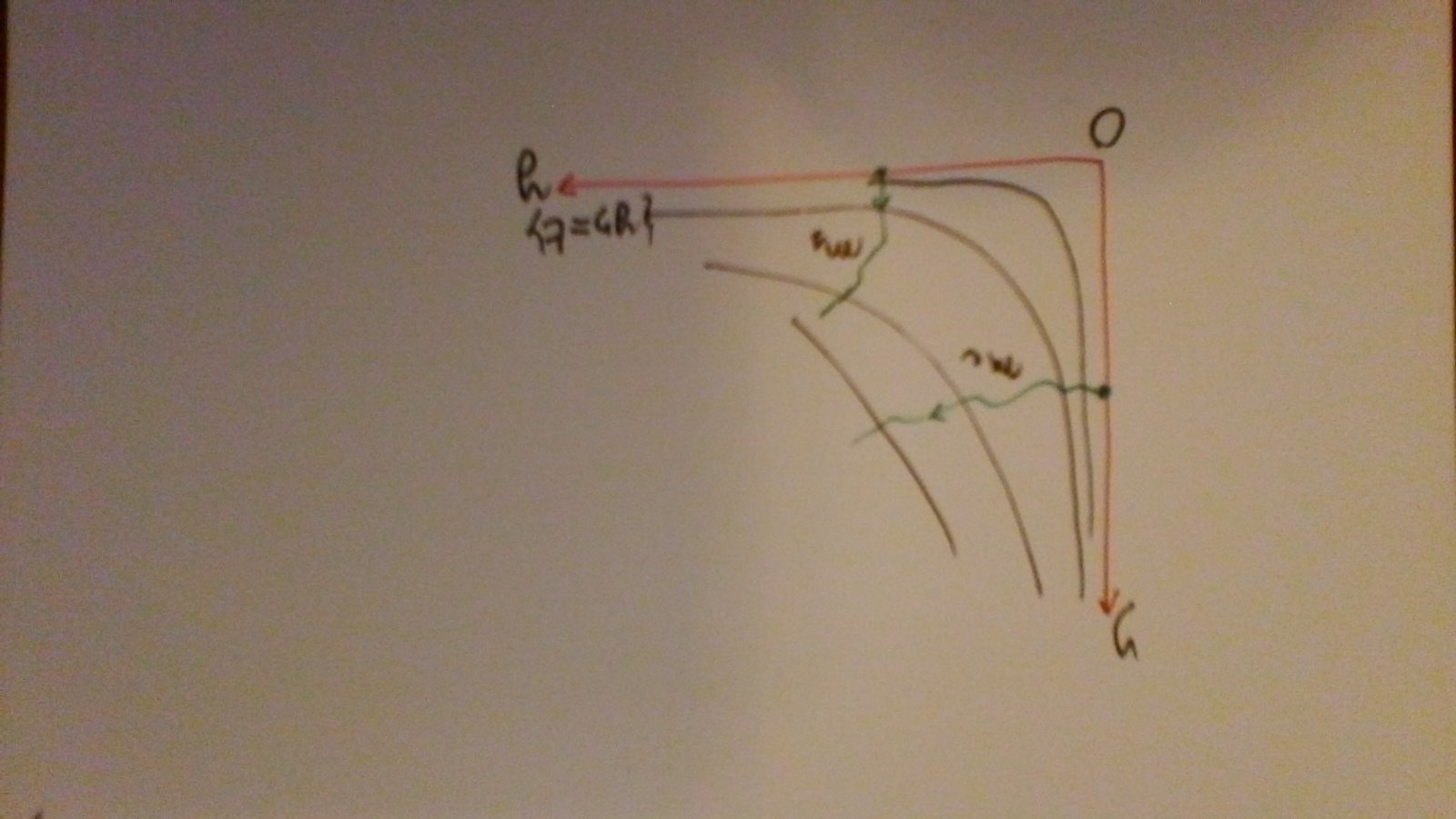}
\caption{the foliation by hyperbolae and the transverse arcs $m_1$ and $m_2$  }
{\label{reduced-1}}
\end{figure}

This lemma follows from the Stokes formula:   
the isochoric Morse lemma (see \cite{Co-Ve-79})  allows to reduce to the case where $d\ga =dy\wedge d\eta $
and to the change of variable  $t\ra F(t)$.
    \subsubsection{The lines $D_j$ as   1D-manifolds with boundary}

Let us put a structure of a 1D-manifold with boundary on  the line $D_j$ in the ``hyperbolic case''.
Let us recall that we denote be $\T_\infty $ the limit of the tori $\T_\gs $ as $\gs $ tends  to one of the infinities.
We showed that   $\T_\infty $ is the union of a finite number of closed hyperbolic orbits 
and of a finite number of cylinders which are parts of the stable an unstable manifolds of these
orbits.
Near  $\T_\infty \setminus Z_0 $, the foliation by the level sets of the moment map is smooth. We can choose any local
transverse arc  to that foliation. They are all equivalent up to diffeomorphism along any connected component of $\T_\infty \setminus Z_0 $
and give local parametrization of $D_j$ near that boundary by intervals $[0,c[$.
How do we pass from one component to the next by  crossing $Z_0$?
  We choose a Poincar\'e section at a point of $Z_0$ and  use the Morse lemma which gives
local coordinates $(y,\eta)$ in that section 
so that $(\gl I +\mu J)(y,\eta) =cte + y\eta $. The local parameter is then the evaluation of the function
$y\eta $ which allows to pass from the transversal $\eta=1$ to the transversal $y=1$. Both are locally parametrized by the restriction
of the function $y\eta$. This gives to $D_j$ the structure of a 1D compact manifold with boundary.
 Note that this holds in a smooth way with respect to  $E'$ close to $E$. 

\subsubsection{The asymptotic behaviour of the action integrals}
There exists, in a  neighborhood $V_j$ of  $\T_{\infty}$, invariant by the flows of $\vec{I}$ and $\vec{J}$,  an  Hamiltonian $P$,
 Poisson
commuting with $I$ and $J$,  
 whose orbits are 
periodic of period $1$ (Theorem 1.6 of \cite{Co-Vu-03}).
$P$ is  constant on $Z_0$.  This   gives  a smooth  action of the group $S^1$ on $V_j$. 
Note that this action is principal on $(V_j \cap \Sigma )\setminus Z_0$, but can get some non trivial isotropy $\Z/2\Z$ on $Z_0$. 
Let    $\gg_1 (z), ~z\in V_j,$ be the 
  $S^1$-orbits.  They are all homotopic. If $z $ lies in some invariant torus, $\gg_1 $ is a homotopically non trivial loop
in this torus.
 We denote by $p_1$ the action integral on $\gg_1(z)$ which is clearly smooth in $V_j$. Note that $p_1$ is a function of $P$ which 
 is  a local diffeomorphism.

We need to choose  a loop $\gg_2$ on the tori in $V_j$   which, with $\gg_1$, generates  a basis of the homology of the invariant tori.
 Let $R_h:=V_j  \cap P^{-1}(h) / S^1 $ with $h$ close to $P(Z_0)$.
The reduced manifolds (see Appendix \ref{app:red}) is foliated by the reduction of the integrable foliation restricted
to $ P^{-1}(h)$. Let us denote by $Z_h$ the singular set of that foliation. 
 As does $V_j\cap \Sigma$, the orbifold 
$R_h$ consists  of a singular part $R_{{\rm sing},h}$, the quotient of $Z_h \cap V_j $, 
 which is homeomorphic to a circle,  and
an open set smoothly foliated by circles which are the reductions of the invariant tori.  Together they  give a topological foliation
of $R_h$ depending smoothly of $h$. 
The singular leaf  $R_{{\rm sing},h}$  is smooth outside the finite set of points which are  quotients by the $S^1$ action 
 of the hyperbolic periodic orbits of $H$.  This foliation is smooth outside these singular points. 
 We take for $\gg_2$ a lift of the projection of $\T_\gs$  depending continuously of $\gs$.

\begin{figure}[!ht]
\centering
\includegraphics[height=7cm]{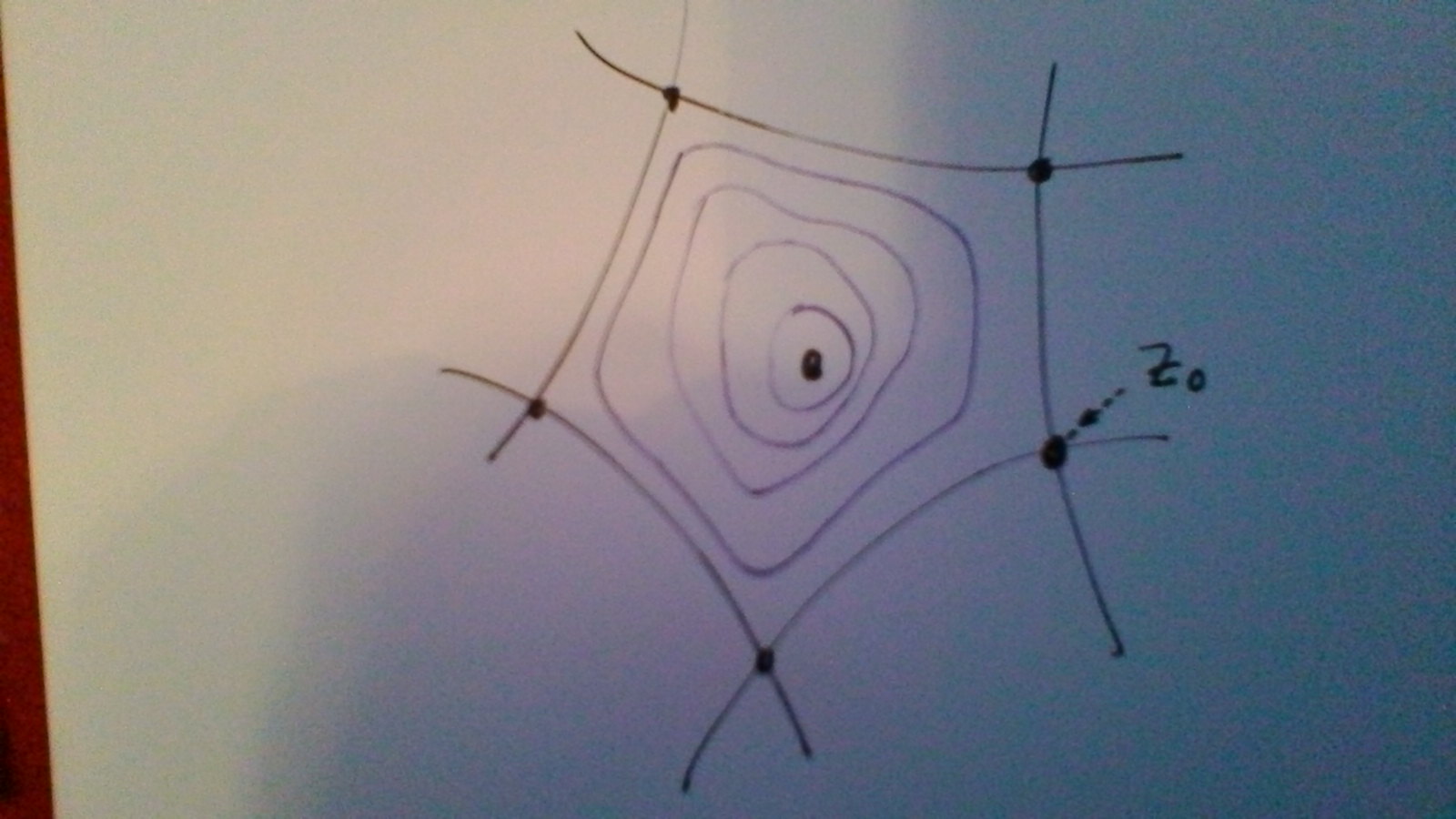}
\caption{a reduced manifold $R_h$}
{\label{reduced}}
\end{figure}

We have the following crucial Lemma:
\begin{lemm} \label{lemm:gradK}
The  action integrals 
 $(p_1,p_2)$ on the previously  choosen  loops $\gg_1$ and $\gg_2$
satisfy at the  boundary
\begin{itemize}
\item 
the action  $p_1$ is smooth up to the  boundary 
\item 
The action $p_2$ as  a function of $\gs $ is of  type {\bf (L)} at the  boundary and  depends smoothly of $h$ and hence of $p_1$. 
\end{itemize}
\end{lemm}
{\it Proof.--}
We saw already the smoothness of $p_1$.
The function  $p_2(\sigma, h) -p_2(Z_h)$ is given by
 the symplectic area in $R_h$ between the reduction of $\gg_2$ and $Z_h$.
 Lemma \ref{lemm:typel} implies that $p_2$ is of type {\bf (L)} depending smoothly of $h$ and hence of $p_1$. 

\hfill  $\square $

\subsubsection{The asymptotic behaviour of $\nabla K $}
We have the following important 
\begin{coro}\label{cor:int-nablaK}
The 1-form  $d\nabla K $ is integrable at any hyperbolic boundary point of ${\cal L}$. 
\end{coro}
{\it Proof.--} 

Near a closed orbit  of $Z_0$, we have the normal form
$H=\Phi(\xi, y\eta )$ with 
$(x,\xi,y,\eta)\in T^\star_{x,\xi} \T \times T^\star_{y,\eta} \R $.
We have $p_1 =\xi $ up to a constant. 
 We get
$H=K(p_1,p_2)=\Phi (p_1, F(p_1, p_2))$ expressing $H$ in terms of the actions.
We get 
\[ \pa_1K =\pa_1 \Phi + \pa_2 \Phi \times  \pa F/\pa p_1,~~
\pa _2 K = \pa_2 \Phi \times  \pa F/\pa p_2 \]
   which are smooth outside $\Sigma $ and continuous on $\Sigma $. Hence their derivatives
are integrable.

\hfill $\square $ 
\section{An integral formula for $|S_t|$}\label{sec:A}
One of the difficulties in extending the result for the disk to this case is the fact that the action-angle coordinates 
only exist outside $Z$. Therefore, we need to make some assumptions on the point $A$. 

\subsection{Assumptions on the point $A$}
If $\Lambda \subset T^\star X $ is a Lagrangian manifold, the {\it caustic set } of $\Lambda $ is the set of critical points
of the projection $\pi_X$ restricted to $\Lambda $. 
We first need a 
\begin{lemm}\label{lemm:caustic}
Let us take  $\go_0\in \Sigma^A$  so that $(A,\go_0)\notin Z_0$ and denote by
 $\F_0$  the 2D-leaf of the invariant foliation of $\Sigma $ containing
$(A,\go_0)$. If $(A,\go_0)$ does not belong to the caustic set of $\F_0$,   then $\Sigma ^A $ and $\F_0$ are transversal
at the point $(A,\go_0)$. 
\end{lemm}
{\it Proof.--} $\Sigma ^A $ is a 1D-submanifold of $T_A^\star X$ and hence $\pi_X (\Sigma ^A)=\{ A \}$. On the other hand, the fact
that $(A,\go_0)$ is not in the caustic set means that 
$(\pi_X)_{|\F_0}$ is a local diffeomorphism onto $X$ near $(A,\go_0)$. The conclusion follows. \hfil $\square$ 

We will  assume: 

{\bf (A3) The intersection  of  $\Sigma ^A$ with $ Z $ is a countable  set.}

\begin{prop}\label{prop:caust1}
 {\bf (A3)} is satisfied as soon as there is only a finite number of $\go \in \Sigma ^A\cap (Z \setminus Z_0)$
 so that $(A,\go)$ is  in the caustic
set of the Lagrangian leaf in which it lives. 
\end{prop}
{\it Proof.--}
The intersection of $\Sigma^A$ with  $Z_0$ is a finite set. On the other hand, the points $\go $ so that 
$(A,\go)$ is not a caustic point of the corresponding leaf are isolated inside $\Sigma^A$. Hence there is at most
a countable set of such points.  \hfil $\square$ 


{\bf (A4) The set  of critical points of the smooth map $\go \ra \gs $ from 
$\Sigma ^A\cap (\Sigma \setminus Z) $ into ${\cal L}$ is countable.}

\begin{prop}
 {\bf (A4)} is satisfied as soon as there is only a finite number of $\go \in \Sigma^A\cap (\Sigma \setminus Z) $
 so that $(A,\go)$ is a caustic point
of the invariant torus  containing that point.  \end{prop}
The argument is quite similar to that of the proof of Proposition \ref{prop:caust1}

\subsection{Exact formulae for $|S_t|$}
We will  compute the lengths of the wave front using the action-angle coordinates. 

We will start with the finite covering of $\Sigma \setminus Z$ by the semi-global  action-angle charts. This allows a description
of $S_t$ as follows: let $\chi :U \ra \T^2 \times C $ be  one of these charts
and let $\Pi_X : \T^2\times C \ra X $ be the map
$\pi _X  \circ \chi^{-1} $.

This way, if we call $(\theta (\go), p( \go ))$ the image of  $\go \in \Sigma ^A$ by $\chi$, we can assume that 
$\theta (\go)$ vanishes identically, because $T^\star_A X$ is Lagrangian.
We get  that the corresponding part of the wave front $S_t$ is defined
by
\[ S_t =\{ \Pi_X \left( t {\rm \tilde{\nabla} }K(p (\go )),p(\go)\right) |\go \in  \Sigma^A\} \] 
where $K$ is the Hamiltonian $H$ expressed in the action coordinates
and $(\theta (\go), p(\go)$ are the action-angle coordinates of $\go \in \Sigma^A$. 
We get, using  the Assumption {\bf (A3)}, the
expression
\[ |S_t|=t\int _{\Sigma ^A} \gamma^\ha \left( \left( {\rm \tilde{\nabla} }K(p (\go )),p(\go)\right); \frac{d}{d\go}
 {\rm \tilde{\nabla} }K(p(\go))\right)d\go 
 \]
where $\gamma $ is the pull-back of $g$ by $\Pi_X$. 
We can make a change of variable: instead of $\go $, one can use $\gs \in {\cal L}$ thanks to assumption  {\bf (A4)}.
We get the
\begin{prop}\label{prop:length} The length of the wave front is given by
\[ |S_t|=t\int _{\cal L}N_A(\gs) \gamma^\ha \left( \left(t {\rm \tilde{\nabla} }K(\sigma ),\gs \right); d \tilde{\nabla }K(\sigma )\right) 
 \]
where $N_A(\gs)=\# \{ \Sigma ^A \cap {\cal L} \}$.
\end{prop}


\section{The main  result}


\begin{theo}\label{theo:main2} Let $(X,g)$ be a   Riemannian manifold of dimension $2$ with $g$ continuous, possibly degenerate. Let 
$H$ be  an  Hamiltonian 
  integrable  at energy $E$ and satisfying the   assumptions {\bf (A1)} and  {\bf(A2)}.
 Let  $A\in X$ be a point 
satisfying the   assumptions {\bf (A3)}  and {\bf   (A4)}.
The length for the metric $g$  of the wave front $S_t $ starting from $A$  has a linear asymptotics
$|S_t|\sim \lambda (A)t $ as $t \ra \infty $.
 
Let us denote be ${\cal L}$ the 1D-manifold of all invariant Lagrangian tori $L_\sigma,~\sigma \in {\cal L},$  filling
  $\Sigma \setminus Z$ and  
 consider the continuous  density on ${\cal L}$ defined
by
\[ |d\sigma |=\int _{\T_\sigma } \gamma^{\ha}\left((\theta,\sigma ); d({\rm \tilde{\nabla}}K (\sigma )) \right)|d\theta |\]
where $\gamma $ is the pull-back of $g$ by the projection $\Pi_X$. 
The measure $|d\sigma |$ is independent of $A$.
 We have
\begin{equation} \label{equ:int2} \lambda (A) =
 \int _{\cal L}N_A(\sigma )|d\sigma | \end{equation}
with $N_A(\sigma ):=\# \{ \Sigma ^A\cap  L_\sigma  \}$. 

\end{theo}
\begin{coro}Let  $ H$ be the Hamiltonian of the
geodesic flow of a smooth metric $G$ on a closed manifold $X$. 
 If $D$ is a smooth domain with boundary in $X$, the $g-$length of $S_t\cap D$, is given
by 
\[ |S_t \cap D|\sim t\int_D d\mu_A  \]
where $d\mu_A $ is an  absolutely continuous density $d\mu_A= F |dx| $ with $F\in L^1 (X,|dx|)$.
 whose integral is $ \lambda (A)$.

\end{coro}
{\it Proof of the Corollary.--}
Let $\psi $ be a positive continuous function on $X$. We can apply the previous theorem with $g'=\psi^2 g$.
This way, we see that the asymptotics of the $g'$-length of $S_t$ is given 
by replacing the measure $|d\sigma |$ by the measure
\[ 
|d\sigma |'=
 \int _{\T_\gs} \psi (\pi_X (\theta, \gs)) \gamma^\ha (\left((\theta,\sigma ); d({\rm \tilde{\nabla}}K (\sigma ))\right) |d\theta |\]
This says that 
\[ |S_t |_{g'}\sim t\int_X \psi d\mu_A \]
where $d\mu_A$ is the pushforward by $\pi_X$ of the absolutely continuous (a.c. in short) finite measure
 $dM_A:=N_A(\gs)\gamma^\ha (\left((\theta,\sigma ); d({\rm \tilde{\nabla}}K (\sigma ))\right) |d\theta | $, supported by $\Sigma $.

We need to show that we can apply this when $\psi $  is the characteristic function of a smooth domain. 
In our situation $\Sigma $ is the unit cotangent bundle and $\pi_X: \Sigma \ra X $ is a submersion. It follows that
 that $d\mu_A $ is a.c. w.r. to $|dx|$. \hfil $\square $
 

\section{Proof of Theorem \ref{theo:main2}}\label{sec:proof}

We start from the expression of $|S_t|$ given in Proposition \ref{prop:length}. 
Let us show that we can apply Lemma \ref{lemm:ergo} to the integral giving $|S_t|/t$. 
In the notations of that lemma, we have
$V(\gs)={\rm \tilde{\nabla} }K(\gs)$.  The Assumption {\bf (A2)} implies that the assumption on $V$ of the Lemma is satisfied.
The function $F$ is given by
$F(\gs, \theta)=N_A(\gs) \gamma^\ha \left( \theta; d \tilde{\nabla }K(\gs)/d\gs)\right) $.
The integrability assumption follows from the 
Corollaries \ref{cor:int-nablaK-ell} and \ref{cor:int-nablaK}
and the upper bound
\[ |N_A(\gs) \gamma^\ha \left( \theta; W \right)|\leq C \| W \| \]
 The continuity with respect to $\theta$ follows from the continuity
of $g$ and the smoothness of  the projection of any $\T_\gs$ onto $X$. It is shown using Lebesgue's dominated convergence Theorem.

\section{Examples}
\subsection{Surfaces of revolution}
Surfaces with a non trivial action of $S^1$ are tori or spheres.
In both case, the metric is given by $g=a(s)^2 d\theta ^2+ ds^2$  where
$s\in \R/L\Z $ in the first case and
$s\in [0,L] $ in the second (in this case $s=0$ and $s=L$ are the poles).

The assumption {\bf (A1)} is satisfied if and only if $a$ is a Morse function.
The assumption {\bf (A2)} is satisfied for a generic $a$.
Assuming {\bf (A1)} and {\bf (A2)}, 
the assumption {\bf (A3)} is satisfied for any point of the torus and for $A$ not a pole in the case of the sphere.
while {\bf (A4)} is always satisfied.
If $A$ is a pole, $|S_t|$ is periodic of period $2L$.

\subsection{Tri-axial ellipso\"ids}
The integrability was found by C. Jacobi (\cite{Ja-39}, see also \cite{Kl-82}
and  section 3.2 of \cite{Co-Vu-03}).
Assumption {\bf (A1)} and {\bf (A2)} are satisfied.
{\bf (A3)} is satisfied for $A$ not an ombilical point while {\bf (A4)} is always satisfied. 
If $A$ is an ombilical point, $|S_t|$ is periodic.




\begin{appendix}

\section{Stationary phase}
For this section, one can look at \cite{Gu-St-77}, chap. 1. 

We want to evaluate the asymptotics as $t\ra +\infty $ of integrals of the
form
\[ I(t):= \int _\R  e^{it S(x)} a(x) dx \]
where $S$ is a real valued smooth function and $a\in C_o^\infty (\R)$. 
We have the
\begin{prop}\label{prop:statnd}
 Let us  assume that the critical points of $S$, i.e. the zeroes of $S'$, are
non degenerate, i.e. $S''(x)\ne 0$ if $S'(x)=0$.
Then, if $x_1,\cdots, x_N$ are the critical points of $S$ in the support of $a$,
 $I(t)$ admits an asymptotic expansion given by 
\begin{equation} \label{equ:stat} I(t) =\sum _{j=1}^{N}
\frac{\sqrt{2\pi}e^{i\ge_j \pi /4}}{|t S'' (x_j)|^{\ha}}e^{itS(x_j)}\left(a(x_j)+ O(t) \right)
  \end{equation}
with $\ge_j = \pm 1$ depending on the sign of $S''_j(0)$. 
\end{prop} 

In the case where the critical points are degenerate, we have the following result:
\begin{prop}\label{prop:statd} If the zeroes of $S'$ in the support of $a$  are of finite order, we have
$I(t) \ra 0$ as $t\ra \infty $.
\end{prop}

Note also that in equation (\ref{equ:stat}), the remainders ``$O(t)$'' are uniform if $S' $ (resp. $a'$) 
 is close to $S$ (resp. close to $a$) in the $C^\infty $ topology and the support of $a$ stays in some fixed bounded intervall.
 

\section{An ``ergodic'' lemma} \label{app:ergo}
This section could be of independent interest.
\begin{lemm}\label{lemm:ergo} 
For $s \in J$ where $J$ is an  interval   of the real line,
 let  $V(s) =A_1(s)\pa _1+A_2(s)\pa _2$ 
 be a  family of constant vector fields on $\T^2$  depending smoothly of $s$. 
Assume that,
for any $s \in J$, there exists two  derivatives $V^{(k)}(s)$ and $V^{(l)}(s)$ which are linearly independent.

Let  $F$ is a function  on  $J \times \T^2$ 
with $F\in C^0(\T^2,L^1(J,ds))$ satisfying the following condition:
there exists a function $\psi \in L^1(J,ds )$ so that
\[ \forall (s,\theta ) \in J\times \T^,~|F(s,\theta) |\leq \psi(s)~.  \]
Then 
\[\lim _{t\ra +\infty} \int _{J} F(s , [tV(s)] ) |ds|  =  \int _{J\times \T^2}F |ds d\theta| \]
\end{lemm}
The assumption on the derivatives of $V$ have the  following geometrical meaning:
 if $V'(s_0)=0$, we get a cusp point which is of finite order;
if $V'(s_0 )\ne 0$, the curvature of the curve $V$ vanishes at a finite order.
In particular the points where $V'$ and $V''$ are linearly dependent are isolated.

{\it Proof.--} It follows from Lebesgue's dominated convergence theorem, that the map 
$f:\theta \ra F(.,\theta) $ is continuous from $\T^2$ into $L^1(J,|ds|)$. 
Let us choose a finite covering of $\T^2$ by balls of centers $\theta _j$, $1\leq j \leq N$,  
so that the $L^1$-oscillation of $f$ is each ball is smaller than $\ge/2$
and a smooth finite partition $(\psi_j)$ of unity subordinated to that covering. 
Let $F_j\in C_o^\infty (J)$ satisfying
$\|F(.,\theta _j )-F_j \|_{L^1}\leq \ge /2$. Such functions do exist (see \cite{Fo-99} Prop. 8.17). 
If  $G(s,\theta ) =\sum _j \psi_j (\theta ) F_j(s) $, we have
\[ \int _J |F(s,\theta )- G(s,\theta )|ds| \leq \ge  \]

This allows to reduce to prove the result for such a function $G$.
We can again approximate $G $ by a function \[L(s,\theta )=\sum _{n\in \Z^2,~\|n\|\leq N } a_n(s){\rm exp}(2\pi i <n| \theta> )\]
uniformly in $L^1 (|ds|)$. We have $a_0(s)=\int _{\T^2}L(s,\theta ) |d\theta| $. 
We are left with the integrals
\[ \int _J a_n(s)e^{2i\pi t<n| V(s)>} |ds| \]
It follows from the assumption on $V$ that such integrals tend to $0$ as $t\ra \infty $ for $n\ne 0$.

\hfill $\square $

\section{Symplectic $S^1$-reduction}\label{app:red}

Let  $P:M \ra \R $ be an Hamiltonian  on a symplectic manifold $(M,\go)$ so that the vector field   $\vec{P}$
is complete and generates an action of $\T$ onto $M$.
Let us assume that this action is almost free: it is free on an open dense subset of $M$ and all the isotropy subgroups are 
finite. Let us look at an  energy shell $S_h:=P^{-1}(h)$ for some $h\in \R$.
The quotient of $S_h $ by the $\T$-action is an orbifold $R_h$. Let us denote by
$\pi_h$ the canonical projection of $S_h$ onto $R_h$.  The orbifold $R_h$ admits an unique symplectic structure $\Omega $ so that
$\pi^\star (\Omega )=\omega $. 

\end{appendix} 

\bibliographystyle{plain}

\end{document}